\RequirePackage{fix-cm}

\documentclass[smallextended,final]{svjour3}       

\usepackage{amssymb,latexsym}
\usepackage{amsmath}

\newtheorem{lema}{Lemma}
\newtheorem{teo}{Theorem}

\journalname{Annali di Matematica}

\begin{document}

\title{A simpler proof of a Katsurada's theorem and rapidly converging series for $\zeta{(2n+1)}$ and $\beta{(2n)}$}

\titlerunning{Rapidly converging series for $\zeta{(2n+1)}$ and $\beta{(2n)}$}

\author{F.~M.~S.~Lima}

\institute{F.M.S.~Lima \at
              Institute of Physics, University of Brasilia, P.O.~Box 04455, 70919-970, Brasilia-DF, Brazil \\
              Tel.: +55-061-31076088\\
              Fax: +55-061-33072363\\
              \email{fabio@fis.unb.br}           
}

\date{Received: 27 June 2013 / Accepted: 28 January 2014}

\maketitle

\begin{abstract}
In a recent work on Euler-type formulae for even Dirichlet beta values, i.e. $\beta{(2n)}$, I have derived an exact closed-form expression for a class of zeta series. From this result, I have conjectured closed-form summations for two families of zeta series.  Here in this work, I begin by using a known formula by Wilton to prove those conjectures. As example of applications, some special cases are explored, yielding rapidly converging series representations for the Ap\'{e}ry constant, $\zeta(3)$, and the Catalan constant, $G = \beta(2)$. Interestingly, our series for $\,\zeta(3)\,$ converges faster than that used by Ap\'{e}ry in his irrationality proof (1978). Also, our series for $\,G\,$ converges faster than a celebrated one discovered by Ramanujan (1915).  At last, I present a simpler, more direct proof for a recent theorem by Katsurada which generalizes the above results.



\end{abstract}

\section{Introduction}

The Riemann zeta function is defined, for $\,\Re{(s)}>1$, as $\,\zeta(s) := \sum_{k=1}^\infty{{\,1/k^s}}$. This function has a singularity (in fact, a simple pole) at $\,s=1$, which corresponds to the divergence of the harmonic series.  For real values of $\,s$, $s>1$, the series converges to a real number between $1$ and $2$, according to the integral test. For positive integer values of $\,s$, $s>1$, one has a well-known formula by Euler (see Ref.~\cite{DeAmo} and references therein):
\begin{equation}
\zeta(2n) = (-1)^{n-1} \, \frac{2^{2n-1} \, B_{2n}}{(2n)!} \; \pi^{2n} ,
\label{eq:Euler}
\end{equation}
$n$ being a positive integer. Here, $\,B_n$ are Bernoulli numbers, i.e.~the rational coefficients of ${\,z^n/n!}\,$ in the Taylor series expansion of ${\,z/(e^z-1)}$, $0 < |z| < 2\,\pi$. For $\,\zeta{(2n+1)}$, on the other hand, no analogous expression is currently known.  This scenario has a `reverse' counterpart on the values of the Dirichlet beta function, defined as $\,\beta(s) := \sum_{k=0}^\infty{(-1)^k/(2\,k+1)^s}$, $s>0$, in the sense that, for integer values of $s$, the following analogue of Eq.~\eqref{eq:Euler} is known:\footnote{Note that this formula remains valid for $\,n=0$, since $\,E_0 =1\,$ and $\,\beta{(1)} = \pi/4$.}
\begin{equation*}
\beta(2n+1) = (-1)^{n} \, \frac{E_{2n}}{2^{2n+2} \, (2n)!} \; \pi^{2n+1} ,
\label{eq:Euler2}
\end{equation*}
where $\,E_{n}$ are Euler numbers, i.e. the (integer) coefficients of ${\,z^n/n!}\,$ in the Taylor expansion of $\, \mathrm{sech}{(z)}$, $|z| < {\,\pi/2}\,$.  For $\,\beta(2n)$, no analogous expression is known, not even for $\,\beta{(2)}$, known as the Catalan's constant $\,G$. Some progress in this direction was reached by K\"{o}lbig (1996), who proved that~\cite{Kolbig}:
\begin{equation}
\beta(2n) = \frac{\psi^{(2n-1)}{\left( \frac14 \right)}}{2\,(2n-1)!\,4^{2n-1}} \, - \frac{(2^{2n}-1) \, |B_{2n}|}{ 2\, (2n)!} \, \pi^{2n},
\label{eq:beta2n}
\end{equation}
where $\psi^{(n)}{(x)}$ is the polygamma function, i.e.~the $n$-th derivative of $\psi{(x)}$, the digamma function.\footnote{The function $\psi{(x)}$, in turn, is defined as the logarithm derivative of $\,\Gamma{(x)}$, the classical gamma function.}
Although this identity resembles Euler's formula, the arithmetic nature of $\,\psi^{(2n-1)}{\left(\frac14 \right)}\,$ is currently unknown.

In a very recent paper, I have succeeded in applying the Dancs and He series expansion method, as introduced in Ref.~\cite{Dancs}, to find similar formulae for $\,\beta(2n)\,$~\cite{Lima2012}. I could then prove that
\begin{eqnarray}
\sum_{k=1}^\infty{\frac{\zeta(2k)}{2k \, (2k+1) \, \ldots \, (2k+2n-1)} \, \left(\frac{1}{2^{2k}}-\frac{1}{4^{2k}}\right)} = \nonumber \\
(-1)^n\,\frac{2^{2n-2}}{\pi^{2n-1}} \, \beta{(2n)} +\frac{n}{(2n)!}\,\ln{2} + \, \frac12 \, \sum_{m=1}^{n-1}{(-1)^m \, \frac{2^{2m}-1}{\pi^{2m}} \,  \frac{\zeta{(2m+1)}}{(2n-2m-1)!}} \, . \quad
\label{eq:my2}
\end{eqnarray}
Since both
\begin{subequations}
\begin{equation}
\sum_{k=1}^\infty{\frac{(2k-1)!}{(2k+2n-1)!} ~ \frac{\zeta(2k)}{2^{\,2k}}}
\label{subeqn1}
\end{equation}
and
\begin{equation}
\sum_{k=1}^\infty{\frac{(2k-1)!}{(2k+2n-1)!} ~ \frac{\zeta(2k)}{4^{\,2k}}}
\label{subeqn2}
\end{equation}
\end{subequations}
converge absolutely, the zeta series in Eq.~\eqref{eq:my2} equals the difference of these individual series. However, the best I could do there in Ref.~\cite{Lima2012} was to investigate the pattern of the analytical results found for the first series, Eq.~\eqref{subeqn1}, for small values of $n$. This did lead me to \emph{conjecture} a formula for its summation which should be valid for any positive integer $\,n$.

Here in this work, I make use of a classical Wilton's formula to prove the above mentioned conjecture. The proof of another conjecture, involving the zeta series in Eq.~\eqref{subeqn2}, follows from Eq.~\eqref{eq:my2}. Finally, on aiming at a generalization of these formulae I substitute the fractions $(1/2)^{2k}$ and $(1/4)^{2k}$ by $\,x^{2k}$, $\,x\,$ being any non-null real number with $\, |x| \le 1$. This has led me to develop a simpler, more direct proof of a theorem by Katsurada (1999)~\cite{Katsurada}.

\section{Zeta series for $\,\zeta(2n+1)\,$ and $\,\beta(2n)$}

The proof of the first conjecture in Ref.~\cite{Lima2012}, involving the zeta series in Eq.~\eqref{subeqn1}, above, follows from a classical result by Wilton for a rapidly convergent series representation for $\,\zeta{(2n+1)}\,$~\cite{Wilton}. Hereafter, we define $\,\,H_n := \sum_{k=1}^n{1/k}\,$ as the $n$-th harmonic number.

\begin{lema}[Wilton's formula]
 \label{le:Wilton}
\;  Let $\,n\,$ be a positive integer.  Then
\begin{eqnarray*}
\frac{\zeta(2n+1)}{(-1)^{\,n-1} \, \pi^{2n}} = \frac{H_{2n+1}-\ln{\pi}}{(2n+1)!} +\sum_{k=1}^{n-1}{\frac{(-1)^k}{(2n-2k+1)!} \, \frac{\zeta{(2k+1)}}{\pi^{2k}}} +\sum_{k=1}^\infty{\frac{(2k-1)!}{(2n+2k+1)!} \, \frac{\zeta(2k)}{2^{\,2k-1}}} .
\end{eqnarray*}

\end{lema}

For a skeleton of the proof, see the original work of Wilton~\cite{Wilton}.  For a more complete proof, see Sec.~4.2 (in particular, pp.~412--413) of Ref.~\cite{LTSri}, a systematic collection of zeta series recently published by Srivastava and Choi.

This lemma allows us to prove the first conjecture of Ref.~\cite{Lima2012}, namely that in its Eq.~(29).

\begin{teo}[First zeta series]
\label{teo:zetaser1}
\;  Let $\,n\,$ be a positive integer. Then
\begin{eqnarray*}
\sum_{k=1}^\infty{\frac{(2k-1)! \; \zeta(2k)}{(2k+2n-1)!} \, \left(\frac{1}{2}\right)^{2k}} = \frac12 \! \left[ \frac{\ln{\pi}-H_{2n-1}}{(2n-1)!} +\sum_{m=1}^{n-1}{(-1)^{m+1} \, \frac{\zeta{(2m+1)}}{\pi^{2m}\,(2n-2m-1)!}} \right] .
\end{eqnarray*}
\end{teo}

\begin{proof}
\;  From Lemma~\ref{le:Wilton}, we know that
\begin{eqnarray*}
\frac{(-1)^n \, \zeta(2n+1)}{\pi^{2n}} = \frac{\ln{\pi} -H_{2n+1}}{(2n+1)!} -\sum_{k=1}^{n-1}{\frac{(-1)^k}{(2n-2k+1)!} \, \frac{\zeta{(2k+1)}}{\pi^{2k}}} \nonumber \\
- \,2 \, \sum_{k=1}^\infty{\frac{(2k-1)!}{(2n+2k+1)!} \, \frac{\zeta(2k)}{2^{\,2k}}} \, .
\end{eqnarray*}
By isolating the last term, one has
\begin{equation*}
2 \sum_{k=1}^\infty{\frac{(2k-1)!}{(2n+2k+1)!} \, \frac{\zeta(2k)}{2^{\,2k}}} = \frac{\ln{\pi} -H_{2n+1}}{(2n+1)!} -\sum_{k=1}^{n}{\frac{(-1)^k}{(2n-2k+1)!} \, \frac{\zeta{(2k+1)}}{\pi^{2k}}} \, .
\end{equation*}
By substituting $\,n = \ell -1\,$ in the above equation, one finds
\begin{equation*}
2 \sum_{k=1}^\infty{\frac{(2k-1)!}{(2k +2\ell-1)!} \, \frac{\zeta(2k)}{2^{\,2k}}} = \frac{\ln{\pi} -H_{2\ell-1}}{(2\ell-1)!} -\sum_{k=1}^{\ell-1}{\frac{(-1)^{k}}{(2\ell-2k-1)!} \, \frac{\zeta{(2k+1)}}{\pi^{2k}}} \, .
\end{equation*}
A division by $2$ completes the proof.
\qed
\end{proof}

For instance, on putting $\,n=2\,$ in this theorem we get a rapidly converging series representation for $\,\zeta(3)$, namely
\begin{equation}
\zeta(3) = \pi^2 \, \sum_{k=1}^\infty{\frac{\zeta(2k)}{k\,(2k + 1)\,(2k + 2)\,(2k + 3) \; 2^{2k}}} +\frac{11}{36}\,\pi^2 -\frac16 \, \pi^2\,\ln{\pi} \, .
\end{equation}
This formula converges to $\zeta(3)$ much faster than $\,\sum_{k\ge1}{1/k^3}$ and even than $\:\frac52 \, \sum_{k \ge 1}{{(-1)^{k+1}} \slash \, {[\,k^3 \, \binom{2k}{k}\,]}}$, a central binomial series used by Ap\'{e}ry~\cite{Apery}. Numerical computation using \emph{Mathematica} shows that only ten terms of the above zeta series are enough for a ten decimal places accuracy.

Now we can use Eq.~\eqref{eq:my2} to prove the other conjecture raised in Ref.~\cite{Lima2012}, namely that in its Eq.~(30).

\begin{teo}[Second zeta series]
\label{teo:zetaser2}
\;  Let $\,n\,$ be a positive integer. Then
\begin{eqnarray*}
\sum_{k=1}^\infty{\frac{(2k-1)! \; \zeta(2k)}{(2k+2n-1)!} \, \left(\frac{1}{4}\right)^{\!2k}} = \frac12 \left[ \frac{\,\ln{(\pi/2)} -H_{2n-1}}{(2n-1)!} -(-1)^{n} \left(\frac{2}{\pi}\right)^{\! 2n-1} \!\! \beta{(2n)} \right. \nonumber \\
\left. \begin{array}{c} ^{} \\ ^{} \end{array}
-\sum_{m=1}^{n-1}{(-1)^m \left(\frac{2}{\pi}\right)^{\!2m} \frac{\zeta{(2m+1)}}{(2n-2m-1)!\,}} \, \right] \! .
\end{eqnarray*}
\end{teo}

\begin{proof}
\;  From Eq.~\eqref{eq:my2}, we know that
\begin{eqnarray*}
\sum_{k=1}^\infty{\frac{(2k-1)! \; \zeta(2k)}{(2k+2n-1)!} \, \left(\frac{1}{2^{2k}}\right)} -\sum_{k=1}^\infty{\frac{(2k-1)! \; \zeta(2k)}{(2k+2n-1)!} \, \left(\frac{1}{4^{2k}}\right)} \nonumber \\
= (-1)^n\,\frac{2^{2n-2}}{\pi^{2n-1}} \, \beta{(2n)} +\frac{n}{(2n)!}\,\ln{2} + \, \frac12 \, \sum_{m=1}^{n-1}{(-1)^m \, \frac{2^{2m}-1}{\pi^{2m}} \,  \frac{\zeta{(2m+1)}}{(2n-2m-1)!}} \: .
\end{eqnarray*}
On substituting the first zeta series at the left-hand side by the result proved in Theorem~\ref{teo:zetaser1}, one has
\begin{eqnarray*}
\frac{\ln{\pi}-H_{2n-1}}{(2n-1)!} -\sum_{m=1}^{n-1}{(-1)^{m} \, \frac{\zeta{(2m+1)}}{\pi^{2m}\,(2n-2m-1)!}} \,-2 \sum_{k=1}^\infty{\frac{(2k-1)! \; \zeta(2k)}{(2k+2n-1)!} \, \left(\frac{1}{4^{2k}}\right)} \nonumber \\
= (-1)^n \left(\frac{2}{\pi}\right)^{2n-1} \! \beta{(2n)} +\frac{\ln{2}}{(2n-1)!} +\, \sum_{m=1}^{n-1}{(-1)^m \, \frac{2^{2m}-1}{\pi^{2m}} \,  \frac{\zeta{(2m+1)}}{(2n-2m-1)!}} \: .
\end{eqnarray*}
By isolating the remaining zeta series, one finds
\begin{eqnarray*}
2 \sum_{k=1}^\infty{\frac{(2k-1)! \: \zeta(2k)}{(2k+2n-1)!} \, \left(\frac{1}{4^{2k}}\right)} = \frac{\ln{(\pi/2)}-H_{2n-1}}{(2n-1)!} -(-1)^n\,\left(\frac{2}{\pi}\right)^{\!2n-1} \! \beta{(2n)} \nonumber \\
-\sum_{m=1}^{n-1}{(-1)^{m} \, \left(\frac{2}{\pi}\right)^{2m} \frac{\zeta{(2m+1)}}{\,(2n-2m-1)! \,}}  \, .
\end{eqnarray*}
A division by $2$ completes the proof.
\qed
\end{proof}

On putting $\,n=1\,$ in Theorem~\ref{teo:zetaser2} we get a rapidly converging series representation for $\,\beta(2) = G$, namely
\begin{equation}
G = \pi \, \sum_{k=1}^\infty{\frac{\zeta(2k)}{2k\,(2k + 1)} \: \frac{1}{4^{2k}}} - \frac{\pi}{2} \: \ln{\!\left(\frac{\pi}{2}\right)}+\frac{\pi}{2} \, .
\end{equation}
This formula converges much faster than $\,\sum_{k \ge 0}{(-1)^k/(2k+1)^2}\,$ and even faster than $\,\frac{\pi}{8} \, \ln{\left( 2 +\sqrt{3} \,\right)} +\frac38 \, \sum_{n \ge 0}{1/\!\left[(2n + 1)^2 \, \binom{2n}{n}\right]}$, a rapidly converging central binomial series discovered by Ramanujan~\cite{rama}. Numerically, only six terms of the above zeta series are enough for a result accurate to ten decimal places. After an extensive search for similar zeta series in literature, I have found a formula by Srivastava and Tsumura (2000) in Ref.~\cite{Sri2000} (see also Ref.~\cite{LTSri}, p.~421, Eq.~(30)\,). In fact, this formula could be taken into account for an independent proof of our Theorem~\ref{teo:zetaser2}, after some simple manipulations, as the reader can easily check.\footnote{\label{ft:polygam} For this, it will be useful to know that $\,\zeta{(s,a)} := \sum_{k=0}^{\infty} 1/(k+a)^s\,$ is the Hurwitz zeta function, for which it is well-known that $\,\zeta{(n+1,a)} = (-1)^{n+1} \, \psi^{(n)}(a)/n!\,$ (see, e.g., Eq.~(25.11.12) of Ref.~\cite{Nist}). Then $\,\zeta{\left( 2m, \frac14 \right)} = \psi^{(2m-1)}\left(\frac14 \right) / (2m-1)!$ and, from our Eq.~\eqref{eq:beta2n}, it follows that $\,\psi^{(2m-1)}{\left( \frac14 \right)} / (2m-1)! = 2^{4m-1} \, \beta(2m) +(-1)^{m-1} \, 2^{4m-2} \, (2^{2m}-1) \, B_{2m} \, \pi^{2m} / (2m)!\,$, where $\,|B_{2m}|\,$ was substituted by $\,(-1)^{m-1} \, B_{2m}$.}

On investigating the substitution of $\,(1/2)^{2k}$ and $\,(1/4)^{2k}$ by $\,x^{2k}$, $\,x$ being any non-null real number with $\,|x| \le 1$, in the zeta series treated in the previous theorems, I have found the following general result.

\begin{teo}[Generalization]
\label{teo:geral}
\quad  Let $\,n\,$ be a positive integer and $\,x\,$ be any real number with $\, 0 < |x| \le 1$. Then
\begin{eqnarray*}
\zeta{(2n+1)} \,- \frac{1}{2\pi x} \, \sum_{\ell=1}^\infty{\frac{\sin{(2 \pi \ell x)}}{\ell^{\,2n+2}}} = (-1)^{n-1} \: (2\pi x)^{2n} \left[ \frac{\,H_{2n+1} -\ln{\! \left( 2 \pi \, |x| \right)}}{(2n+1)!} \right. \nonumber \\
\left. \begin{array}{c} ^{} \\ ^{} \end{array}
+ \sum_{k=1}^{n-1}{(-1)^k \frac{\zeta{(2k+1)}}{(2n-2k+1)! \; (2\pi x)^{2k}}} \; + 2 \sum_{k=1}^\infty{\frac{(2k-1)! \: \zeta{(2k)}}{(2n+2k+1)!} ~ x^{\,2k}\,} \right] \! .
\end{eqnarray*}
\end{teo}

\begin{proof}
From the Euler's product formula for the sine function, we know that, for all non-null real $\,z\,$ with $\,|z| \le 1$, the following identity holds:
\begin{equation}
\frac{\sin{(\pi z)}}{\pi z} = \prod_{n=1}^\infty{\left( 1-\frac{z^2}{n^2} \right)} .
\end{equation}
On taking the logarithm on both sides, we have
\begin{equation}
\ln{\left[\frac{\sin{(\pi z)}}{\pi z}\right]} = \sum_{n=1}^\infty{\ln{\left( 1-\frac{z^2}{n^2} \right)}} = - \, \sum_{n=1}^\infty{ \left( \sum_{k=1}^\infty{ \frac{z^{2k}}{k \: n^{2k}} }\right)} ,
\end{equation}
where the logarithm was expanded in a Taylor series in the last step. This implies that
\begin{equation}
\ln{\left[\frac{\,2 \, \sin{(\pi z)}}{2 \, \pi z}\right]} = - \, \sum_{k=1}^\infty{ \left( \, \sum_{n=1}^\infty{ \frac{1}{n^{2k}} } \right) \, \frac{z^{2k}}{k}} = -\sum_{k=1}^\infty{\frac{z^{2k}}{k} \, \, \zeta{(2k)}}
\end{equation}
and then
\begin{equation}
\ln{\left| 2 \, \sin{(\pi z)}\right|} - \ln{\left|2 \, \pi z\right|} = \, - \, \sum_{k=1}^\infty{\frac{z^{2k}}{k} \, \, \zeta{(2k)}} \, .
\end{equation}
On multiplying both sides by $\,(x-z)^{2n}$, $n$ being a positive integer, and integrating from $0$ to $x$, $x$ being any non-null real in the interval $[-1,1]$, we have
\begin{eqnarray}
\int_0^x{(x-z)^{2n} \, \ln{\left| 2 \, \sin{(\pi z)}\right|} \: dz} \, - \int_0^x{(x-z)^{2n} \, \ln{\left| 2 \, \pi z\right|} \: dz} \nonumber \\
= \, - \, \int_0^x{(x-z)^{2n} \, \sum_{k=1}^\infty{\frac{z^{2k}}{k} \, \zeta{(2k)}} \: dz} \, .
\label{eq:Zints}
\end{eqnarray}
Let us solve each definite integral carefully.

The first integral in Eq.~\eqref{eq:Zints} can be expanded in a trigonometric series as follows. Since, for all $\,\theta \in (0,2\pi)$,\footnote{This is a well-known Fourier series expansion.}
\begin{equation}
\ln{\left[2 \, \sin{\left(\frac{\theta}{2}\right)} \right]} = - \sum_{k=1}^\infty{\frac{\cos{(k \, \theta)}}{k}} \, ,
\end{equation}
then, on substituting $\,z = \theta / (2\pi)\,$ in the integral, one finds
\begin{eqnarray}
I_n(x) := \int_0^x{(x-z)^{2n} \: \ln{\left| 2 \, \sin{(\pi z)}\right|} \, dz} \nonumber \\
= \frac{1}{(2\pi)^{2n+1}} \, \int_0^{2 \pi x}{(2 \pi x -\theta)^{2n} \: \ln{\left| 2 \, \sin{\left(\frac{\theta}{2}\right)} \right|} \: d \theta } \nonumber \\
= \, - \, \frac{1}{(2\pi)^{2n+1}} \, \int_0^{2 \pi x}{(2 \pi x -\theta)^{2n} \: \sum_{k=1}^\infty{\frac{\cos{(k \, \theta)}}{k}} \: d \theta } \nonumber \\
= \, - \, \frac{1}{(2\pi)^{2n+1}} \, \int_0^{2 \pi x}{(2 \pi x -\theta)^{2n} \: d\left(\sum_{k=1}^\infty{\frac{\sin{(k \, \theta)}}{k^2}}\right) } . \:
\end{eqnarray}
On integrating by parts, one has
\begin{eqnarray}
I_n(x) = \frac{2n}{(2\pi)^{2n+1}} \, \int_0^{2 \pi x}{(2 \pi x -\theta)^{2n-1} \; d\left(\sum_{k=1}^\infty{\frac{\cos{(k \, \theta)}}{k^3}}\right) } \nonumber \\
= - \frac{2n}{(2 \pi)^2} \, x^{2n-1} \, \zeta(3) + \frac{(2n)!}{(2\pi)^{2n+1} \, (2n-2)!} \, \int_0^{2 \pi x}{(2 \pi x -\theta)^{2n-2} \; d\left(\sum_{\,k=1}^\infty{\frac{\sin{(k \, \theta)}}{k^4}}\right) } . \qquad {}
\end{eqnarray}
On integrating by parts again and again, we find, after some algebra, that
\begin{eqnarray}
I_n(x) = \sum_{j=1}^n{\frac{(-1)^j \, (2n)! \, \zeta(2j+1)}{(2\pi)^{2j} \, (2n+1-2j)!} \, x^{2n+1-2j} } \, + \, \frac{(-1)^{n-1} \, (2n)!}{(2\pi)^{2n+1}} \, \sum_{k=1}^\infty{\frac{\sin{(2 k \, \pi \, x)}}{k^{2n+2}}} \, . \qquad
\label{eq:Zint1}
\end{eqnarray}

The second integral in Eq.~\eqref{eq:Zints} is readily solved by noting that
\begin{eqnarray}
\int_0^x{(x-z)^{2n} \, \ln{| 2 \, \pi z |} \: dz} = \, - \, \frac{1}{2n+1} \, \int_0^x{\ln{|2\pi\,z|} \, \, d \left[ (x-z)^{2n+1} -x^{2n+1} \right]}  \nonumber \\
= \, \frac{x^{2n+1}}{2n+1} \, \left[ \ln{|2\,\pi\,x|} + \sum_{\ell=1}^{2n+1}{\frac{(-1)^\ell}{\ell} \: \binom{2n+1}{\ell}} \right] . \quad
\label{eq:Zint2}
\end{eqnarray}

The third integral, i.e.~the one at the right-hand side of Eq.~\eqref{eq:Zints}, can be written in the form of a zeta series on integrating it by parts directly. This yields
\begin{equation}
\int_0^x{(x-z)^{2n} \, \sum_{k=1}^\infty{\frac{z^{2k}}{k} \, \zeta{(2k)}} \: dz}  = 2 \, (2n)! \, \sum_{k=1}^\infty{\frac{(2k-1)! \, \zeta(2k) \, x^{2k+2n+1}}{(2k+2n+1)!}} \, .
\label{eq:Zint3}
\end{equation}

Finally, on substituting the results in Eqs.~\eqref{eq:Zint1}, \eqref{eq:Zint2}, and \eqref{eq:Zint3} on the integrals in Eq.~\eqref{eq:Zints}, we find
\begin{eqnarray}
\sum_{j=1}^n{\frac{(-1)^j \, (2n)! \, \zeta(2j+1)}{(2\pi)^{2j} \, (2n+1-2j)!} \, \frac{1}{x^{2j}} } \, + \, \frac{(-1)^{n-1} \, (2n)!}{(2\pi \, x)^{2n+1}} \, \sum_{k=1}^\infty{\frac{\sin{(2 k \, \pi \, x)}}{k^{2n+2}}}  \nonumber \\
-\frac{1}{2n+1} \left[ \ln{|2\,\pi\,x|} + \sum_{\ell=1}^{2n+1}{\frac{(-1)^\ell}{\ell} \: \binom{2n+1}{\ell}} \right]  =  - 2 \, (2n)! \, \sum_{k=1}^\infty{\frac{(2k-1)! \, \zeta(2k)}{(2k+2n+1)!} \: x^{2k}} . \qquad
\label{eq:Zexpr}
\end{eqnarray}
On dividing both sides by $[-(2n)!]$, we have
\begin{eqnarray}
2 \sum_{k=1}^\infty{\frac{(2k-1)!}{(2k+2n+1)!} \, \zeta(2k) \: x^{2k}} = \, - \, \sum_{j=1}^n{\frac{(-1)^j \, \zeta(2j+1)}{(2\pi \, x)^{2j} \, (2n+1-2j)!} } \, + \frac{\ln{(2\,\pi\,|x|)} - H_{2n+1}}{(2n+1)!} \nonumber \\
+ \, \frac{(-1)^n}{(2 \pi \, x)^{2n+1}} \, \sum_{k=1}^\infty{\frac{\sin{(2 k \, \pi \, x)}}{k^{2n+2}}} \, , \qquad {}
\label{eq:Zexprfim}
\end{eqnarray}
where we have made use of the binomial sum $\,H_{2n+1} = \, - \, \sum_{\ell=1}^{2n+1}{\frac{(-1)^\ell}{\ell} \: \binom{2n+1}{\ell}}$, which is easily proved by induction on $\,n$. On extracting the last term of the sum involving odd zeta values (i.e., that for $j=n$) we bring up the desired result.
\qed
\end{proof}

Our proof of Theorem~\ref{teo:geral}, above, has allowed us to detect some typos in Theorem~2 of Ref.~\cite{Katsurada}, as well as some mistakes in its proof, which is based upon the Mellin transform technique.  In fact, the formula as printed in Katsurada's paper cannot be correct because as $\,x \rightarrow 0\,$ the left-hand side approaches $\,2 \, \zeta(2n+1)\,$ whereas the right-hand side approaches zero. Unfortunately, that incorrect formula is reproduced on p.~442 of Ref.~\cite{LTSri}, which is currently the main reference on zeta series in literature. Our Theorem~\ref{teo:geral} corrects both the absence of a modulus in the argument of the logarithm and the `$+$' sign preceding the sine series in Eq.~(1.6) of Ref.~\cite{Katsurada}.

Interestingly, the formula in our Theorem~\ref{teo:geral} can be reduced to a more appropriate form for both symbolic and numerical computations. This is easily obtained by noting that $\,\sum_{\ell=1}^\infty{\sin{(2 \pi x \, \ell)}/\ell^{\,2n}} = \mathrm{Cl}_{2n}(2\pi x)\,$ for all real $x$, where $\,\mathrm{Cl}_{2n}(\theta) := \Im{\left\{\mathrm{Li}_{2n}\left(e^{i\,\theta}\right)\right\}} \,$ is the Clausen function of order $\,2 n$. By substituting $\:n = m-1\:$ and $\:x = \theta/(2 \pi)\,$ in Theorem~\ref{teo:geral}, we readily find  \begin{eqnarray}
\sum_{k=1}^\infty{\frac{(2k-1)! ~ \zeta{(2k)}}{(2m+2k-1)!} \;  {\left( \frac{\theta}{2 \pi}\right)\!}^{2k}} = \, \frac12 \left[\frac{\:\ln{|\theta|} -H_{2m-1}}{(2m-1)!} \,-(-1)^m \, \frac{\mathrm{Cl}_{2m}(\theta)}{\theta^{2m-1}} \right. \nonumber \\
\left. \begin{array}{c} ^{} \\ ^{} \end{array}
 -\sum_{k=1}^{m-1}{(-1)^k \, \frac{\zeta{(2k+1)}}{(2m-2k-1)! \, ~ \theta^{2k}}} \, \right] \! ,
\label{eq:Cl2m}
\end{eqnarray}
which holds for any non-null real $\,\theta\,$ with $\, |\theta| \le 2\pi$. Advantageously, this form remains valid for all \emph{positive integer} values of $\,m$, as long as the sum at the right-hand side is taken as null when $\,m=1$, as usual. Since $\,\mathrm{Cl}_{2m}(\pi) = 0\,$ and $\,\mathrm{Cl}_{2m}(\pi/2) = \beta{(2m)}$, Eq.~\eqref{eq:Cl2m} allows for prompt proofs of Theorems~\ref{teo:zetaser1} and~\ref{teo:zetaser2}, respectively, showing that they are special cases of Theorem~\ref{teo:geral}, as expected. This is indeed the case for a number of rapidly convergent zeta series in literature, as e.g.~some of the zeta series given in Refs.~\cite{Borwein,Dancs,Lima2012,Sri98,Wilton} and many zeta series presented in Chaps.~3 and 4 of Ref.~\cite{LTSri}. This reflects the generality of our Theorem~\ref{teo:geral}. Equation~\eqref{eq:Cl2m} can thus be viewed as a source of rapidly converging zeta series for odd zeta values and Clausen functions of even order.

For instance, on taking $\,\theta = \pi/3\,$ in Eq.~\eqref{eq:Cl2m}, one finds
\begin{eqnarray}
2 \, \sum_{k=1}^\infty{\frac{(2k-1)!}{(2m+2k-1)!} \;  \frac{\,\zeta{(2k)}}{\,6^{2k}}} = \frac{\:\ln{(\pi/3)} -H_{2m-1}}{(2m-1)!} \,-(-1)^m \, \frac{\,\mathrm{Cl}_{2m}(\pi/3)}{\,(\pi/3)^{2m-1}} \nonumber \\
-\sum_{k=1}^{m-1}{(-1)^k \, \frac{\zeta{(2k+1)}}{(2m-2k-1)! \, ~ (\pi/3)^{2k}}} \, . \quad
\end{eqnarray}
For $\,m = 1$, this reduces to
\begin{equation}
\sum_{k=1}^\infty{\frac{\zeta{(2k)}}{k \, (2k+1)} \;  \frac{1}{6^{2k}} } = \ln{\!\left(\frac{\pi}{3}\right)} -1 +\frac{3}{\pi} \: \mathrm{Cl}_{2}\!\left(\frac{\pi}{3}\right) ,
\label{eq:pi3m1}
\end{equation}
whereas for $\,m = 2\,$ one has
\begin{eqnarray}
\sum_{k=1}^\infty{\frac{\zeta{(2k)}}{k\,(2k+1)\,(2k+2)\,(2k+3)} \;  \frac{1}{\,6^{2k}}} = \frac16 \, \ln{\!\left(\frac{\pi}{3}\right)} -\frac{11}{36} \,- \frac{27}{\,\pi^3} \: \mathrm{Cl}_{4}\!\left(\frac{\pi}{3} \right) + \frac{9}{\,\pi^2} \; \zeta{(3)} \, . \qquad
\label{eq:pi3m2}
\end{eqnarray}
Similarly, for $\,\theta = \pi/4\,$ one has
\begin{equation}
\sum_{k=1}^\infty{\frac{\zeta{(2k)}}{k \, (2k+1)} \;  \frac{1}{8^{2k}} } = \ln{\!\left(\frac{\pi}{4}\right)} -1 +\frac{4}{\,\pi} \; \mathrm{Cl}_{2}\!\left(\frac{\pi}{4}\right)
\label{eq:pi4m1}
\end{equation}
for $\,m = 1\,$ and
\begin{equation}
\sum_{k=1}^\infty{\frac{\zeta{(2k)}}{k\,(2k+1)\,(2k+2)\,(2k+3)} \;  \frac{1}{8^{2k}} } = \frac16 \, \ln{\!\left(\frac{\pi}{4} \right)} -\frac{11}{36} -\frac{64}{\,\pi^3} \: \mathrm{Cl}_{4}\!\left( \frac{\pi}{4} \right)   +\frac{16}{\,\pi^2} \; \zeta{(3)}
\label{eq:pi4m2}
\end{equation}
for $\,m = 2$. I could not find these formulae explicitly in literature.

\section{Rates of convergence}
All zeta series investigated here belong to the class embraced by Theorem~\ref{teo:geral}.  Their rate of convergence can be analyzed as follows. For convenience, denote by $\,\mathcal{S}_k\,$ the summand of the zeta series at the right-hand side of the formula in our Theorem~\ref{teo:geral}. By applying Stirling's formula, namely $\,k! \sim \left( k/e \right)^k \sqrt{2 \pi k}\,$, and noting that $\,\zeta{(2k)} \rightarrow 1\,$ as $\,k \rightarrow \infty$, one finds
\begin{eqnarray}
\mathcal{S}_k \sim \frac{(2k-1)^{2k-1}}{(2n+2k+1)^{2n+2k+1}}  ~ e^{2n+2} \, \frac{\sqrt{2k-1}}{\sqrt{2n+2k+1}} ~ x^{2k} \nonumber \\
= \frac{1}{\left(1+\frac{2n+2}{2k-1}\right)^{2k-1}}  ~ e^{2n+2} \, \frac{1}{\sqrt{1+\frac{2n+2}{2k-1}}} \, \frac{1}{(2k+2n+1)^{2n+2}} ~ x^{2k} \nonumber \\
\sim \frac{1}{e^{2n+2}} ~ e^{2n+2} \, \frac{1}{\sqrt{1+\frac{2n+2}{2k-1}}} \, \frac{1}{(2k+2n+1)^{2n+2}} ~ x^{2k} \nonumber \\
= \frac{1}{\sqrt{1+\frac{2n+2}{2k-1}}} \, \frac{1}{(2k+2n+1)^{2n+2}} ~ x^{2k} \nonumber \\
\sim \frac{1}{(2k+2n+1)^{2n+2}} ~ x^{2k} \sim \frac{1}{(2k)^{2n+2}} ~ x^{2k} \, ,
\end{eqnarray}
where we have made use of $\,\lim_{\,y \rightarrow \infty} \left(1+\alpha/y \right)^y = e^\alpha\,$ and the binomial approximation $\,(1+y)^\alpha \approx 1+\alpha\,y\,$. Therefore,
\begin{equation}
\mathcal{S}_k = O \! \left(\frac{x^{2k}}{(2k)^{2n+2}} \right) \qquad (k \rightarrow \infty, \, n \in \mathbb{N}) \, ,
\end{equation}
which is valid for any real $\,x\,$ with $\,|x| \le 1$.  This allows for a direct comparison of rates of convergence for the zeta series investigated here in this paper.  For instance, the zeta series in Wilton's formula has $\, x = \frac12 \,$, hence $\,\mathcal{S}_k = O \! \left({2^{-2k} \: (2k)^{-2n-2}} \right) = O \! \left({2^{-2k-2n-2} \: k^{-2n-2}} \right) = O \! \left({2^{-2k} \: k^{-2n-2}} \right)$, whereas the zeta series in our Theorem~\ref{teo:zetaser2} has $\,x=\frac14\,$, hence $\,\mathcal{S}_k = O \! \left({4^{-2k} \: (2k)^{-2n-2}} \right) = O \! \left({4^{-2k-n-1} \: k^{-2n-2}} \right) = O \! \left({4^{-2k} \: k^{-2n-2}} \right)$. This is why the latter converges somewhat faster than the former.

\begin{acknowledgements}
The author would like to thank Mrs.~Claudionara de Carvalho for many interesting and useful discussions.
\end{acknowledgements}


\begin{thebibliography}{}

\bibitem{Apery} Ap\'{e}ry R.: Irrationalit\'{e} de $\zeta(2)$ et $\,\zeta(3)$. Ast\'{e}risque \textbf{61}, 11--13 (1979)

\bibitem{Borwein} Borwein, J.M., Bradley, D.M., Crandall, R.E.: Computational strategies for the Riemann zeta function. J. Computat. Appl. Math. \textbf{121}, 247--296 (2000)

\bibitem{Dancs} Dancs, M.J., He, T.-X.: An Euler-type formula for $\zeta{(2k+1)}$. J. Number Theory \textbf{118}, 192--199 (2006)

\bibitem{DeAmo} De Amo, E., Diaz Carrillo, M., Fernandez-Sanchez, J.: Another proof of Euler's formula for $\zeta(2k)$. Proc. Amer. Math. Soc. \textbf{139}, 1441--1444 (2011)

\bibitem{Katsurada} Katsurada, M.: Rapidly convergent series representations for $\zeta(2n+1)$ and their $\chi$-analogue. Acta Arithmetica \textbf{90}, 79--89 (1999)

\bibitem{Kolbig} K\"{o}lbig, K.S.: The polygamma function $\psi^{(k)}(x)$ for $x =1/4$ and $x =3/4$. J. Comput. Appl. Math. \textbf{75}, 43--46 (1996)

\bibitem{Lima2012} Lima, F.M.S.: An Euler-type formula for $\,\beta(2n)\,$ and closed-form expressions for a class of zeta series. Int. Transf. Special Funct. \textbf{23}, 649--657 (2012)

\bibitem{Nist} Olver, F.W.J., Lozier, D.W., Boisvert, R.F., Clark, C.W.: \textsl{NIST Handbook of Mathematical Functions}. Cambridge University Press, New York (2010)

\bibitem{rama} Ramanujan, S.: On the integral $\int_0^x{\frac{\tan^{-1}{t}}{t} \, dt}$. J. Indian Math. Soc. \textbf{VII}, 93--96 (1915)

\bibitem{Sri98} Srivastava, H.M.: Further series representations for some definite integrals associated with the Riemann zeta function. Appl. Math. Comput. \textbf{97}, 1--15 (1998)

\bibitem{Sri2000} Srivastava, H.M., Tsumura, H.: A certain class of rapidly convergent series representations for $\zeta{(2n+1)}$. J. Comput. Appl. Math. \textbf{118}, 323--335 (2000)

\bibitem{LTSri} Srivastava, H.M., Choi, J.: \textsl{Zeta and $q$-zeta functions and associated series and integrals}. Elsevier, Amsterdam (2012)

\bibitem{Wilton} Wilton, J.R.: A proof of Burnside's formula for $\log{\Gamma(x+1)}$ and certain allied properties of Riemann's $\zeta$-function. Messenger Math. \textbf{52}, 90--93 (1922/1923)

\end{thebibliography}
\end{document}